\documentclass[12pt]{amsart}
\usepackage{amsmath}
\usepackage{amscd}
\usepackage{graphics}
\usepackage{latexsym}

\theoremstyle{plain}

\newtheorem{Thm}{Theorem}
\newtheorem{Lem}[Thm]{Lemma}
\newtheorem{Cor}[Thm]{Corollary}
\newtheorem{Prop}[Thm]{Proposition}
\newtheorem{Def}[Thm]{Definition}

\newcommand{\OB}{\mathfrak{{ob}}}

\newcommand{\bfr}{{\mathbb{R}}}

\newcommand{\Li}{{\mathbb{L}}}

    \def\sqr#1#2{{\vcenter{\hrule height.#2pt
            \hbox{\vrule width.#2pt height#1pt \kern#1pt
            \vrule width.#2pt}\hrule height.#2pt}}}
    \def\square{\mathchoice\sqr67\sqr67\sqr{2.1}6\sqr{1.5}6}

\begin{document}

\title[]{An open book decomposition compatible with rational contact surgery}
\author{Burak Ozbagci}
\address{Department of Mathematics \\ Ko\c{c} University \\ Istanbul, Turkey}
\email{bozbagci@ku.edu.tr}

\begin{abstract}
We construct an open book decomposition compatible with a contact
structure given by a rational contact surgery on a Legendrian link
in the standard contact $S^3$. As an application we show that some
rational contact surgeries on certain Legendrian knots induce
overtwisted contact structures.

\end{abstract}
\maketitle
\setcounter{section}{-1}
\section{Introduction}

\noindent Recently Giroux proved a central result regarding the
topology of contact 3-manifolds. Namely he established a one to
one correspondence between contact structures up to isotopy and
open book decompositions up to positive stabilizations. This
correspondence, however, does not explicitly describe an open book
decomposition corresponding to a given contact structure.

In \cite{dg}, Ding and Geiges proved that every (closed) contact
3-manifold $(Y, \xi)$ can be given by a contact $(\pm 1)$-surgery
on a Legendrian link in the standard contact $S^3$. Here we use
the parenthesis to emphasize that the surgery coefficients are
measured with respect to the contact framing. If the coefficients
of all the curves in a contact surgery diagram are $(-1)$, then an
open book decomposition compatible with this contact structure is
given by the algorithm in \cite{ao} coupled with the work of
Plamenevskaya (\cite{pl}). Moreover, Stipsicz (\cite{sti}) showed
that the same algorithm works in the general case of contact $(\pm
1)$-surgery. In this article we will review these results (giving
slightly different proofs) and extend the algorithm to the case of
rational contact surgery.

In fact any rational contact surgery can be turned into a sequence
of contact $(\pm 1)$-surgeries (\cite{dg, dgs}) and the algorithm
above would provide an open book decomposition compatible with the
resulting contact structure. However, this would give an open book
decomposition with high genus and we will show that there is a
short-cut in obtaining an open book decomposition (with lower
genus) compatible with a rational contact surgery.

As an application we show that certain rational contact surgeries
induce overtwisted contact structures by making use of the
right-veering property of tight contact structures recently
introduced by Honda, Kazez and Matic (\cite{hkm}).

Here we outline the main idea in the article: Given a Legendrian
link in the standard contact $(S^3, \xi_{st})$ with its front
projection onto the $yz$-plane. We will show that there is an open
book decomposition on $S^3$ compatible with a contact structure
$\xi_0$ isotopic to $\xi_{st}$ such that the Legendrian link
(after a Legendrian isotopy) is contained in a page of this open
book decomposition as described in \cite{pl}. It follows that the
page framing and the contact framing coincide on each component of
this link since the Reeb vector field induces both framings.
Consequently when we perform contact $(\pm 1)$-surgery on this
link we get an open book decomposition compatible with the
resulting contact structure. The monodromy of this open book
decomposition is given by a product of Dehn twists along curves
explicitly drawn on a page. Now when a Legendrian link with
rational surgery coefficients is given,  we embed this link into
the page of the open book decomposition as described above. Then
we attach appropriate 1-handles to a page of this open book
decomposition and extend the monodromy of our open book
decomposition by Dehn twists along some push-offs of the original
monodromy curves going through the attached 1-handles.

\;

\noindent {\bf {Acknowledgement}}: We would like to thank Olga
Plamenevskaya, Andras Stipsicz and John Etnyre for helpful
correspondence. The author was partially supported by the Turkish
Academy of Sciences.

\section{Open book decompositions and contact structures}

We will assume that all our contact structures are positive and
cooriented. In the following we describe the compatibility of an
open book decomposition with a given contact structure on a
3-manifold.

{\Def \label{openbook}
  Suppose that for a link $L$ in a $3$-manifold $Y$ the complement $Y-L$
  fibers as $\pi \colon Y-L \to S^1$ such that the fibers are interiors of
  Seifert surfaces of $L$. Then $(L, \pi)$ is an \emph{open book
    decomposition} of $Y$. The Seifert surface $F={\overline {\pi^{-1}(t)}}$ is called
    a \emph{page}, while $L$
the \emph{binding} of the open book decomposition. The monodromy
of the fibration $\pi$ is called the \emph{monodromy} of the open
book decomposition.}

\;

Any locally trivial bundle with fiber $F$ over an oriented circle
is canonically isomorphic to the fibration $I\times F/(1,x)\sim
\big (0,h(x)\big ) \to I/\partial I\approx S^1$ for some
self-diffeomorphism $h$ of $F$. In fact, the map $h$ is determined
by the fibration up to isotopy and conjugation by an orientation
preserving self-diffeomorphism of $F$. The isotopy class
represented by the map $h$ is called the monodromy of the
fibration. Conversely given a compact oriented surface $F$ with
nonempty boundary and $h \in \Gamma_F$ (the mapping class group of
$F$) we can form the mapping torus $F(h)=I\times F/(1,x)\sim \big
(0, h(x)\big )$. Since $h$ is the identity on $\partial F$, the
boundary $\partial F(h)$ of the mapping torus $F(h)$ can be
canonically identified with $r$ copies of $T^2 = S^1 \times S^1$,
where the first $S^1$ factor is identified with $I/\partial I$ and
the second one comes from a component of $\partial F$. Hence by
gluing in $r$ copies of $D^2\times S^1$ to $F(h)$ so that
$\partial D^2$ is identified with $S^1 = I/\partial I$ and the
$S^1$ factor in $D^2 \times S^1$ is identified with a boundary
component of $\partial F$, $F(h)$ can be completed to a closed
$3$-manifold $Y$ equipped with an open book decomposition. In
conclusion, an element $h\in \Gamma_F$ determines a $3$-manifold
together with an ``abstract" open book decomposition on it.

{\Thm [Alexander] Every closed and oriented 3-manifold admits an
open book decomposition.}

Suppose that an open book decomposition with page $F$ is specified
by $h\in \Gamma_F$. Attach a $1$-handle to the surface $F$
connecting two points on $\partial F$ to obtain a new surface
$F^{\prime}$. Let $\alpha$ be a closed curve in $F^{\prime}$ going
over the new $1$-handle exactly once. Define a new open book
decomposition with $ h \circ t_\alpha \in \Gamma_{F^{\prime}} $,
where $t_\alpha$ denotes the right-handed Dehn twist along
$\alpha$. The resulting open book decomposition is called a
\emph{positive stabilization} of the one defined by $h$. If we use
a left-handed Dehn twist instead then we call the result a
\emph{negative stabilization}. The inverse of the above process is
called \emph{positive} (\emph{negative}) \emph{destabilization}.
Note that the topology of the underlying 3-manifold does not
change when we stabilize/destabilize an open book. Also note that
the resulting monodromy depends on the chosen curve $\alpha$.

{\Def \label{compatible} An open book decomposition of a
$3$-manifold $Y$ and a contact structure $\xi$ on $Y$ are called
\emph{compatible} if $\xi$ can be represented by a contact form
$\alpha$ such that the binding is a transverse link, $d \alpha$ is
a volume form on every page and the orientation of the transverse
binding induced by $\alpha$ agrees with the boundary orientation
of the pages.}

\;

In other words, the conditions that $\alpha>0 $ on the binding and
$d\alpha>0 $ on the pages is a strengthening of the contact
condition $\alpha \wedge d \alpha >0$ in the presence of an open
book decomposition on $Y$. The condition that $d \alpha$ is a
volume form on every page is equivalent to the condition that the
Reeb vector field of $\alpha$ is transverse to the pages. Moreover
an open book decomposition and a contact structure are compatible
if and only if the Reeb vector field of $\alpha$ is transverse to
the pages (in their interiors) and tangent to the binding.

{\Thm [Giroux] Every contact 3-manifold admits a compatible open
book decomposition.}

\section{An open book decomposition compatible with a contact $(\pm 1)$-surgery}
\label{comp}

In this section we describe an explicit construction of an open
book decomposition compatible with a given contact structure. The
algorithm is contained in \cite{ao} and it is proven to be
compatible with the given contact structure in \cite{pl} and
\cite{sti}.

We will show that for a given Legendrian link $\Li$ in $(\bfr^3,
\xi_{st}) \subset (S^3, \xi_st)$ there exists a surface $F \subset
S^3$ containing $\Li$ such that $d\alpha$ is an area form on $F$
(where $\alpha= dz +x\, dy$), $\partial F =K$ is a torus knot
which is transverse to $\xi_{st} = \ker d\alpha$ and the
components of $\Li$ do not separate $F$. We first isotope $\Li$ by
a Legendrian isotopy so that in the front projection (onto the
$yz$-plane) all the segments have slope $(\pm 1)$ away from the
points where $\Li$ intersects the $yz$-plane. Then we consider
narrow rectangular strips around each of these segments and
connect them by small twisted bands corresponding to each point
where $\Li$ intersects the $yz$-plane. The small bands can be
constructed in such a way that the Legendrian link lies on these
bands while the bands twist along the contact planes. The narrow
strips around the straight segments connected with these small
twisted bands give us the Seifert surface $F$ of a torus knot
$K=\partial F$. Notice that we ensured that $\Li$ lies in $F$.
Moreover $d\alpha$ is an area form on $F$ by construction since
the Reeb vector field $\frac{\partial}{\partial z}$ of $\alpha$ is
transverse to $F$. Furthermore we can slightly isotope $\partial
F=K$ to make it transverse to $\xi_{st}$.

Now since $K$ is a fibered knot with fibered surface $F$ there is
a fibration of the complement of $K$ in $S^3$ where $F$ is one of
the pages of the induced open book decomposition on $S^3$. Note
that $d\alpha$ induces an area form on the nearby pages as well
since we can always keep the nearby pages transverse to
$\frac{\partial}{\partial z}$. The union of these nearby pages
including the binding $K$ is a handlebody $U_1$ (which is a
thickening of this one page $F$ that carries the Legendrian link
$\Li$) such that $d\alpha$ is an area form on every page. But we
can not guarantee that $d\alpha$ induces an area form on the rest
of the pages of this open book decomposition. We would like to
extend the contact structure $\xi_{st}$ to the complementary
handlebody $U_2$ (as some contact structure $\xi_0$) so that it is
compatible with the pages in $U_2$. This can be achieved (see
\cite{pl}) by an explicit construction of a contact form on $U_2$
similar to the one described in \cite{tw}.  Hence we get a contact
structure $\xi_0$ on $S^3$ which is compatible with our open book.
Moreover, by construction $\xi_0$ and $\xi_{st}$ coincide on $U_1$
and we claim that the contact structures $\xi_0$ and $ \xi_{st}$
are isotopic on $U_2$ relative to $\partial U_2 $. Notice that
$\partial U_2 $ can be made convex and one can check that the
binding $K$ is the dividing set on $\partial U_2 $. Uniqueness (up
to isotopy) of a tight contact structure with such boundary
conditions is given by Theorem~\ref{torisu}.





Suppose that $(K,\pi) $ is a given open book decomposition on a
closed $3$-manifold $Y$. Then by presenting the circle $S^1$ as
the union of two closed (connected) arcs $S^1=I_1\cup I_2$
intersecting each other in two points, the open book decomposition
$(K,\pi)$ naturally induces a \emph{Heegaard decomposition}
$Y=U_1\cup_{\Sigma } U_2$ of the $3$-manifold $Y$. The surface
$\Sigma$ along which these handlebodies are glued is simply the
union of two pages $\pi^{-1}(I_1\cap I_2)$ together with the
binding.

{\Thm [Torisu, \cite{to}] \label{torisu} Suppose that $\xi_1$,
$\xi_2$ are contact structures on $Y$ satisfying:

\noindent \emph{(i)} ${\xi_i|}_{U_j}$ $(i=1,2;$ $j=1,2)$ are
tight, and

\noindent \emph{(ii)} $\Sigma$ is convex in $(Y, \xi_i)$ and $K$
is the dividing set for both contact structures.

\noindent Then $\xi_1$ and $\xi_2$ are isotopic. In addition, the
set of such contact structures is nonempty.}

\;

Summarizing the above discussion we get

{\Prop  [Plamenevskaya, \cite{pl}] \label{olgap} For a given
Legendrian link $\Li$ in $(S^3 , \xi_{st})$ there exists an open
book decomposition of $S^3$ satisfying the following conditions:

\noindent $(1)$ the contact structure $\xi_0$ compatible with this
open book decomposition is isotopic to $\xi_{st}$,

\noindent $(2)$ $\Li$ is contained in one of the pages and none of
the components of\/ $\Li$ separate $F$,

\noindent $(3)$ $\Li$ is Legendrian with respect to $\xi_0$,

\noindent $(4)$ there is an isotopy which fixes $\Li$ and takes
$\xi_0$ to $\xi_{st}$,

\noindent $(5)$ the page framing of\/ $\Li$ (induced by $F$) is
the same as its contact framing induced by $\xi_0$ (or
$\xi_{st}$).}

\;

\noindent In fact item $(5)$ in the theorem above follows from
$(1)$-$(4)$ by

{\Lem \label{agree} Let $C$ be a Legendrian curve on a page of a
compatible open book decomposition $\OB_\xi$ in a contact
3-manifold $(Y,\xi)$. Then the page framing of $C$ is the same as
its contact framing.}

\begin{proof}
Let $\alpha$ be the contact 1-form for $\xi$ such that $\alpha >
0$ on the binding and $d\alpha >0 $ on the pages of $\OB_\xi$.
Then the Reeb vector field $R_\alpha$ is transverse to the pages
(in their interiors) as well as to the contact planes. Hence
$R_\alpha$ defines both the page framing and the contact framing
on $C$.
\end{proof}

Given a Legendrian link $\Li$ in $(\bfr^3, \xi_{st}) \subset (S^3,
\xi_{st})$ we described an open book decomposition on $S^3$ whose
page is the Seifert surface of an appropriate torus knot $K$ and
$\Li$ is included in one of the pages. When we perform contact
$(\pm 1)$-surgery along $\Li$ we get a new open book decomposition
on the resulting contact 3-manifold obtained by the surgery. The
monodromy of this open book decomposition is given by the
composition of the monodromy of the torus knot and Dehn twists
along the components of the surgery link. Here all the Dehn twists
of the monodromy of the torus knot is right-handed while a
$(+1)$-surgery curve (resp. $(-1)$) induces a left-handed (resp.
right-handed) Dehn twist. Notice that the surgery curves are
pairwise disjoint and they are homologically non-trivial on the
Seifert surface by this construction. It turns out that the
resulting contact 3-manifold and the open book decomposition are
compatible by the following theorem a proof of which is can be
found in \cite{gay} in case of $(-1)$-surgery and \cite{et1} for
the general case.

{\Prop \label{surg} Let $C$ be a Legendrian curve on a page of a
compatible open book decomposition $\OB_\xi$ with monodromy $h \in
\Gamma_F$ on a contact 3-manifold $(Y,\xi)$. Then the contact
3-manifold obtained by contact $(\pm 1)$-surgery along $C$ is
compatible with the open book decomposition with monodromy $h
\circ (t_C)^{\mp 1} \in \Gamma_F $, where $t_C$ denotes a
right-handed Dehn twist along $C$.}

\section{An open book decomposition compatible with a rational contact surgery}
\label{r}

 In this section we will first outline how to turn a
rational contact surgery into a sequence of contact $(\pm
1)$-surgeries. The reader is advised to turn to \cite{dg, dgs} for
background on contact surgery.

Assume that we want to perform contact $(r)$-surgery on a
Legendrian knot $L$ in $(S^3, \xi_{st})$ for some rational number
$r<0$. In this case the surgery can be replaced by a sequence of
contact $(-1)$-surgeries along Legendrian knots associated to $L$
as follows: suppose that $r=-\frac{p}{q}$ and the continued
fraction coefficients of $-\frac{p}{q}$ are equal to $[r_0+1, r_1,
\ldots , r_k]$, with $r_i \leq -2$ ($i=0, \ldots , k $). Consider
a Legendrian push-off of $L$, add $|r_0 +2|$ zig-zags to it and
get $L_0$. Push this knot off along the contact framing and add
$|r_1+2|$ zig-zags to it to get $L_1$. Perform contact
$(-1)$-surgery on $L_0$ and repeat the process with $L_1$. After
$(k+1)$ steps we end up with a diagram involving only contact
$(-1)$-surgeries. The result of the sequence of contact
$(-1)$-surgeries is the same as the result of the original contact
$(r)$-surgery according to \cite{dg, dgs}.

{\Prop [Ding--Geiges, \cite{dg}] \label{k} Fix $r=\frac{p}{q}>0 $
and an integer $k>0$. Then contact $(r)$-surgery on the Legendrian
knot $L$ is the same as contact $(\frac{1}{k})$-surgery on $L$
followed by contact $(\frac{p}{q-kp})$-surgery on the Legendrian
push-off\/ $L'$ of\/~$L$.}

\;

By choosing $k>0$ large enough, the above proposition provides a
way to reduce a contact $(r)$-surgery (with $r>0$) to a contact
$(\frac{1}{k})$-surgery and a negative contact $(r')$-surgery.
This latter one can be turned into a sequence of contact
$(-1)$-surgeries, hence the algorithm is complete once we know how
to turn contact $(\frac{1}{k})$-surgery into contact $(\pm
1)$-surgeries.

{\Lem [Ding--Geiges, \cite{dg}] Let $L_1, \cdots , L_k$ denote $k$
Legendrian push-offs of the Legendrian knot $L$. Then contact
$(\frac{1}{k})$-surgery on $L$ is isotopic to performing contact
$(+1)$-surgeries on the $k$ Legendrian knots $L_1, \cdots , L_k$.}

\;

Given a Legendrian link $\Li$ in $(\bfr^3, \xi_{st}) \subset (S^3,
\xi_{st})$ with rational surgery coefficients. We follow the
algorithm described in Section~\ref{comp} to find an open book
decomposition on $S^3$ such that $\Li$ is embedded in one of the
pages. Now use the above algorithm to turn the rational surgery
into contact $(\pm 1)$-surgeries. Since the contact framing of
each component of $\Li$ agrees with the page framing by
Lemma~\ref{agree}, a contact push-off of any component will still
lie on the same page.

\begin{figure}[ht]
  \begin{center}
     \includegraphics{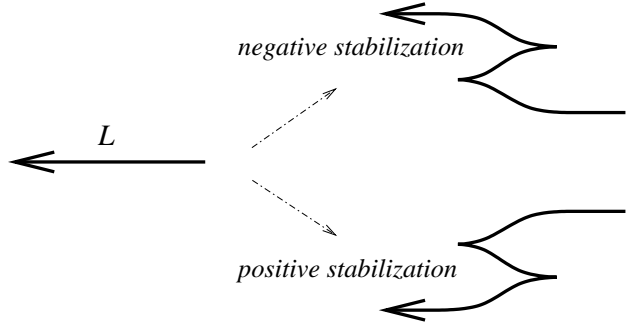}
   \caption{Positive and negative stabilization of a Legendrian knot $L$}
   \label{zig}
    \end{center}
  \end{figure}

Let $L$ be a Legendrian knot in $(\mathbb{R}^3, \xi_{st})$. We
define the positive and negative stabilization of $L$ as follows:
First we orient the knot $L$ and then if we replace a strand of
the knot by an up (down, resp.) cusp by adding a zigzag as in
Figure~\ref{zig} we call the resulting Legendrian knot the
negative (positive, resp.) stabilization of $L$. Notice that
stabilization is a well defined operation, i.e., it does depend at
what point the stabilization is done.

\begin{figure}[ht]
  \begin{center}
     \includegraphics{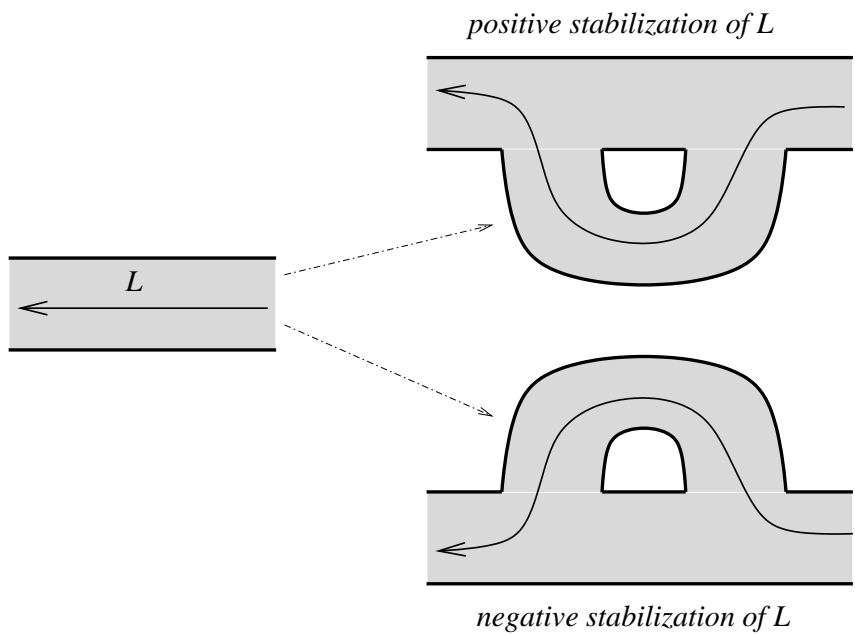}
   \caption{Stabilization of a page to include the stabilization of the
   Legendrian knot $L$} \label{legenstab}
    \end{center}
  \end{figure}

Now let $L$ be a Legendrian knot in a page of a compatible open
book decomposition $\OB_\xi$ in a contact 3-manifold $(Y, \xi)$.
Then by Lemma 3.3 in \cite{et2}, the stabilized knot lies in a
page of an open book decomposition obtained by stabilizing
$\OB_\xi$ by attaching a 1-handle to the page of $\OB_\xi$ and
letting the stabilized knot go through the 1-handle once. Notice
that there is a positive and a negative stabilization of the
oriented Legendrian knot $L$ defined by adding a down or an up
cusp, and this choice corresponds to adding a  \emph{left} (i.e.,
to the left-hand side of the oriented curve $L$) or a \emph{right}
(i.e., to the right-hand side of the oriented curve $L$) 1-handle
to the surface respectively as shown in Figure~\ref{legenstab}.

Repeating this process (by attaching appropriate right or left
1-handles to the Seifert surface of the torus knot) we will get a
page of an open book decomposition where we have all the push-offs
with their additional zig-zags embedded in this page. Notice that
contact $(r)$-surgery is not uniquely defined because of the
choice of adding up or down zig-zags to the push-offs and this
choice can be followed in the way that we attach our 1-handles as
in Figure~\ref{legenstab}. The monodromy of the resulting open
book decomposition will be the composition of the monodromy of the
torus knot (a product of right-handed Dehn twists), right-handed
Dehn twists corresponding to the stabilizations and the Dehn
twists along the push-offs. This open book decomposition is
compatible with the contact $(r)$-surgery by
Proposition~\ref{surg} since we can recover the affect of this
rational surgery on $(S^3, \xi_{st})$ by contact $(\pm
1)$-surgeries along embedded curves on a page of an open book
decomposition of $S^3$ compatible with its standard contact
structure. Here notice that when we positively stabilize a
compatible open book decomposition of $(S^3, \xi_{st})$, the
resulting open book decomposition (of $S^3$) will be compatible
with $\xi_{st}$. As a result we get

{\Thm Given a contact 3-manifold obtained by a rational contact
surgery on a Legendrian link in the standard contact $S^3$. Then
there is an algorithm to find an open book decomposition on this
3-manifold compatible with the contact structure.}

\section{An example}

\begin{figure}[ht]
  \begin{center}
     \includegraphics{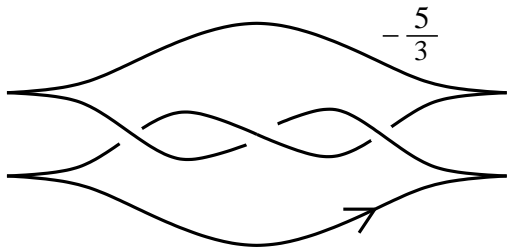}
   \caption{Rational contact surgery on a right-handed Legendrian trefoil knot in $(S^3, \xi_{st})$} \label{rat}
    \end{center}
  \end{figure}

Next we illustrate the algorithm on a simple example. Consider a
contact $(-\frac{5}{3})$-surgery on the  right-handed Legendrian
trefoil knot $L$ in Figure~\ref{rat}.

\begin{figure}[ht]
  \begin{center}
     \includegraphics{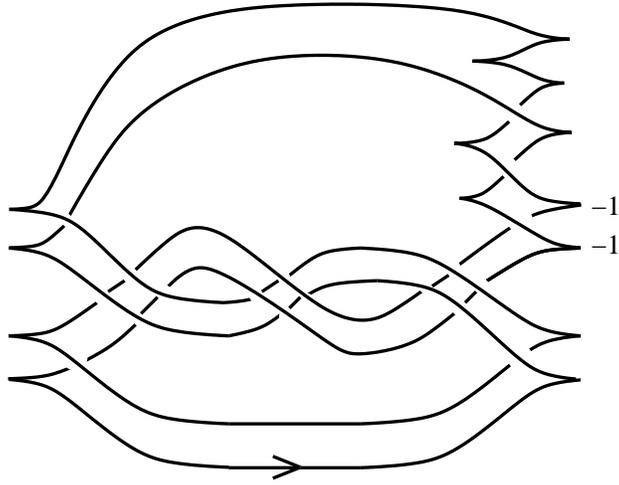}
   \caption{Turning rational contact surgery into contact $(\pm 1)$-surgeries} \label{push}
    \end{center}
  \end{figure}

Observe first that the continued fraction coefficients of
$-\frac{5}{3}$ are $r_0 = r_1 = -3$.  Notice that this surgery is
not uniquely defined and there are four different possibilities of
performing this surgery. We will find an open book decomposition
compatible with one of these surgeries. First orient the
Legendrian knot as indicated in Figure~\ref{rat}. Consider a
Legendrian push-off of $L$, add a down zig-zag to it and get
$L_0$. Push $L_0$ off along the contact framing and add a down
zig-zag to it to get $L_1$. Performing contact $(-1)$-surgery on
both $L_0$ and $L_1$ is equivalent to performing a contact
$(-\frac{5}{3})$-surgery on the right-handed Legendrian trefoil
knot in Figure~\ref{rat}. The push-offs $L_0$ and $L_1$ are
illustrated in Figure~\ref{push}.

\begin{figure}[ht]
  \begin{center}
     \includegraphics{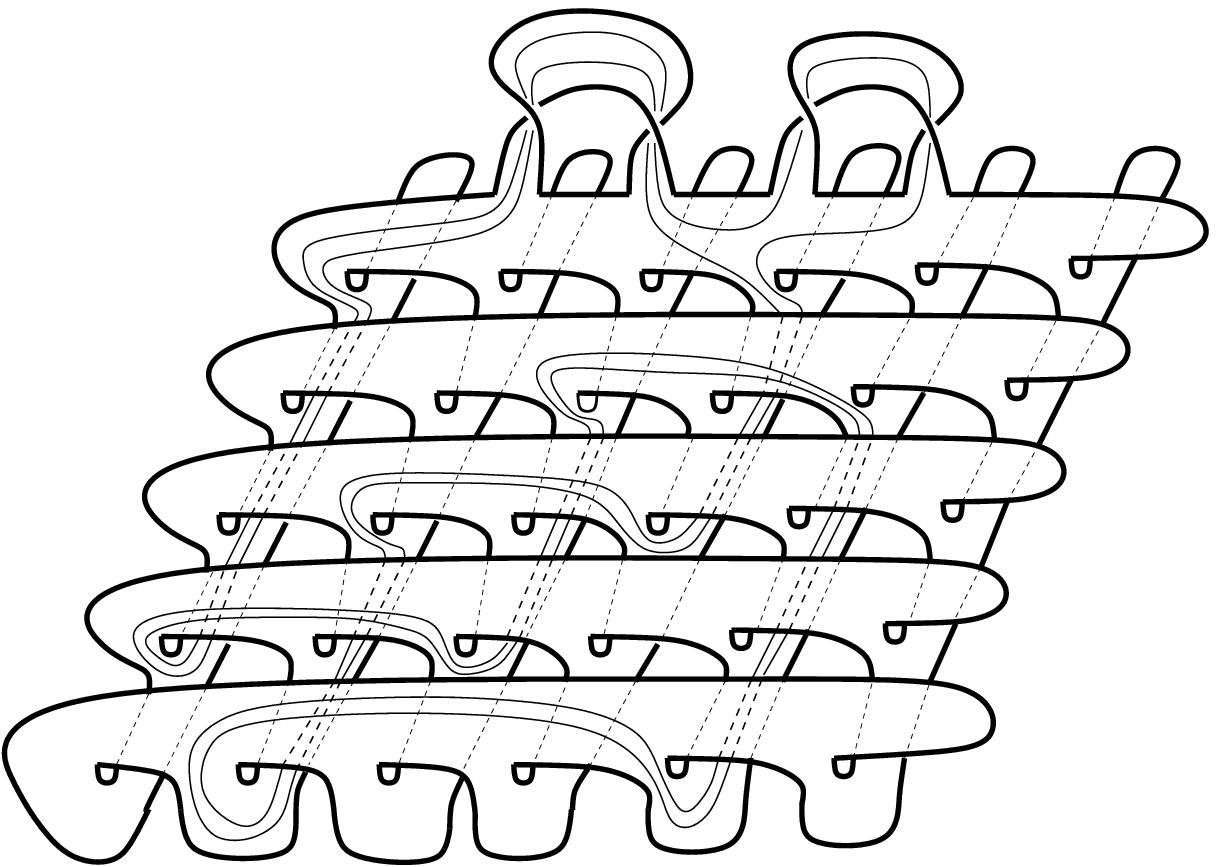}
   \caption{} \label{56}
    \end{center}
  \end{figure}

Finally in Figure~\ref{56}, we depict the page of the open book
decomposition which is compatible with the contact
$(-\frac{5}{3})$-surgery on the  right-handed Legendrian trefoil
knot $L$. The page is the Seifert surface of the (5,6)-torus knot
with two 1-handles attached. Notice that the 1-handle attachments
in Figure~\ref{legenstab} are shown abstractly but in
Figure~\ref{56} this corresponds to plumbing positive Hopf-bands
to the Seifert surface which is embedded in $\mathbb{R}^3$. The
push-offs $L_0$ and $L_1$ are embedded on this page and the
monodromy of the open book decomposition is the product of the
monodromy of the (5,6)-torus knot (see \cite{o} for details),
right-handed Dehn twists along the embedded curves $L_0$ and $L_1$
and right-handed Dehn twists along the core circles of the
1-handles.

\section{An application}

In \cite{o}, we proved that for any positive integer $k$ contact
$(\frac{1}{k})$-surgery on a stabilized Legendrian knot in the
standard contact $S^3$ induces an overtwisted contact structure,
using sobering arcs introduced by Goodman \cite{good}. Note that
contact $(\frac{1}{k})$-surgery is uniquely defined for any
integer $k \neq 0$. In this section we will prove that for any
positive rational number $r$, at least one of the contact
$(r)$-surgeries on a stabilized Legendrian knot induces an
overtwisted contact structure, using the the following

{\Thm [Honda--Kazez--Matic, \cite{hkm}] \label{hkm} If a contact
3-manifold $(Y, \xi)$ is tight then every open book of $Y$
compatible with $\xi$ is right-veering.}

\;

\begin{figure}[ht]
  \begin{center}
     \includegraphics{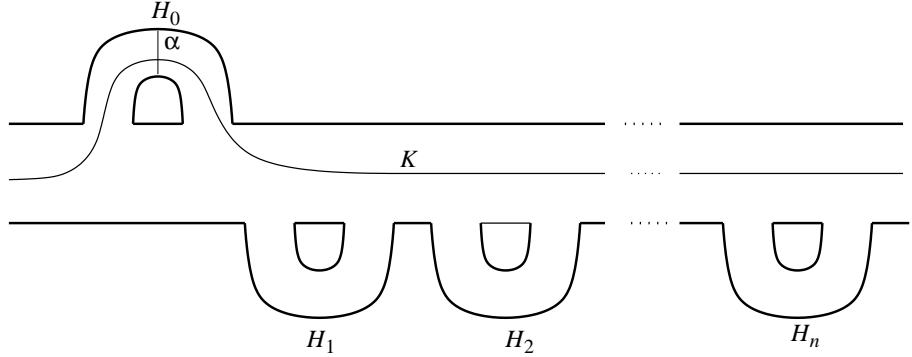}
   \caption{Legendrian knot $K$ on the page and the arc $\alpha$ across $H_0$} \label{handles}
    \end{center}
  \end{figure}

\noindent Suppose that $K$ is a positively stabilized Legendrian
knot in $(S^3, \xi_{st})$. Let $r$ be a positive rational number
and apply the algorithm in Section~\ref{r} to turn a contact
$(r)$-surgery on $K$ into a sequence of contact $(\pm
1)$-surgeries along some push-offs of $K$ in such a way that all
the push-offs are only negatively stabilized. Now consider the
open book we described in Section~\ref{r} compatible with this
surgery diagram. Since $K$ is already positively stabilized, the
page of the compatible open book is obtained by attaching a left
1-handle $H_0$ and some right 1-handles $H_1,H_2, \cdots, H_n$
(corresponding to the push-offs of $K$) to the Seifert surface of
an appropriate torus knot. Note that all the surgery curves are
embedded disjointly on the resulting page of the open book. We
claim that the arc $\alpha$ across $H_0$ (depicted in
Figure~\ref{handles}) is not right-veering at its top-end and
hence showing that the induced contact structure is overtwisted by
the criterion given in Theorem~\ref{hkm}.

To verify that $\alpha$ is not right-veering at its top-end we
apply the monodromy $h$ of the open book to $ \alpha$ and observe
that $h(\alpha)$ lies to the left of $\alpha$ on the page. Here
note that the relevant part of the monodromy $h$ consists of a
product of $k$ (described by Proposition~\ref{k}) left-handed Dehn
twists along $K$, right-handed Dehn twists along the push-offs of
$K$ (going through the 1-handles in Figure~\ref{handles} in
various ways) and a right-handed Dehn twist along the core circle
of the handle $H_0$. The key point is that all the curves along
which we apply right-handed Dehn twists stay only on one side of
$K$ on the surface and a cancellation of a left-handed Dehn twist
by a right-handed Dehn twist is not allowed by construction. Here
we would like to point out that $\alpha$ is not a sobering arc, so
we could not use Goodman's criterion (cf. \cite{good}) to prove
overtwistness of the induced contact structure as we did in
\cite{o}.

\;

Thus we proved

{\Prop \label{stote} For any positive rational number $r$, at
least one of the contact $(r)$-surgeries on a stabilized
Legendrian knot in the standard contact 3-sphere induces an
overtwisted contact structure.}

\;

\begin{figure}[ht]
  \begin{center}
     \includegraphics{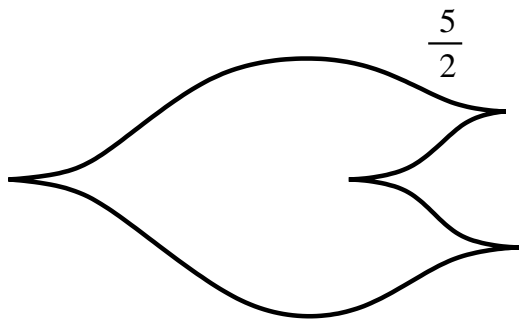}
   \caption{$(\frac{5}{2})$-contact surgery on a Legendrian unknot} \label{fivetwo}
    \end{center}
  \end{figure}

We depict in Figure~\ref{stot} an example of a contact structure
which is overtwisted by Proposition~\ref{stote}. The contact
structure in Figure~\ref{stot} is obtained by a
$(\frac{5}{2})$-contact surgery on the Legendrian unknot which is
shown in Figure~\ref{fivetwo}. The next result easily follows from
Proposition~\ref{stote}:

\begin{figure}[ht]
  \begin{center}
     \includegraphics{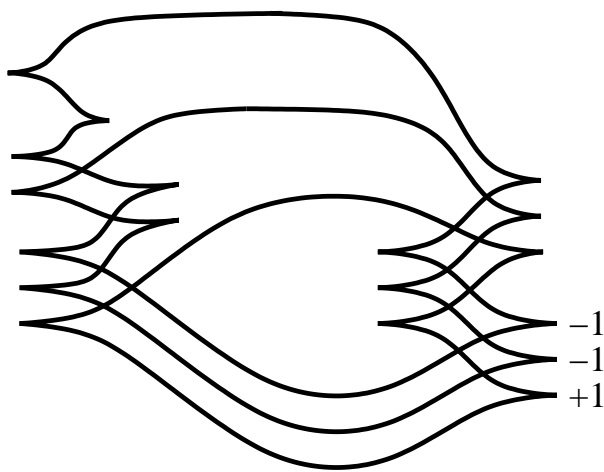}
   \caption{An overtwisted contact structure} \label{stot}
    \end{center}
  \end{figure}

{\Cor If a rational contact surgery diagram contains a Legendrian
knot with an isolated stabilized arc whose surgery coefficient is
positive then at least one of the surgeries it represents is
overtwisted.}

\end{document}